\documentclass[10pt]{amsart}
\usepackage{amsmath,amssymb,graphicx}
\usepackage[all]{xy}
\usepackage[bbgreekl]{mathbbol}
\usepackage{verbatim}
\usepackage{enumerate}

\newtheorem{thm}{Theorem}[section]

\theoremstyle{definition}

\numberwithin{equation}{section}

\DeclareFontFamily{U}{rsf}{} \DeclareFontShape{U}{rsf}{m}{n}{
  <5> <6> rsfs5 <7> <8> <9> rsfs7 <10->  rsfs10}{}
\DeclareMathAlphabet{\mathscr}{U}{rsf}{m}{n}

\renewcommand{\imath}{\sqrt{-1}}

\newcommand{\LL}{\mathbb{L}}

\newcommand{\TT}{\mathbb{T}}

\numberwithin{equation}{section}

\begin{document}

\title[]{Multiple Derived Lagrangian Intersections}
\author{Oren Ben-Bassat}\thanks{}%
\address{Department of Mathematics,
University of Haifa,
Haifa, Israel}
\email{ben-bassat@math.haifa.ac.il}
\dedicatory{}
\thanks{I would like to thank Dominic Joyce for his conjecture on the existence of this Lagrangian structure and many key conversations. I would also like to thank Lino Amorim, Chris Brav, Dennis Borisov, Vittoria Bussi, Damien Calaque, Kobi Kremnizer and Tony Pantev for helpful discussions. I acknowledge the support of the European Commission under the Marie Curie Programme for the IEF grant which enabled this research to take place. The contents of this article reflect my own views and not the views of the European Commission.}

\begin{abstract}
We give a new way to produce examples of Lagrangians in shifted symplectic derived stacks, based on multiple intersections.  Specifically, we show that an $m$-fold  homotopy fiber product of Lagrangians in a shifted symplectic derived stack its itself Lagrangian in a certain cyclic product of pairwise homotopy fiber products of the Lagrangians. 
\noindent 
\end{abstract}
\maketitle 
\section{Introduction}
In a recent article \cite{PTVV} of Pantev, To\"{e}n, Vaqui\'{e}, and Vezzosi, shifted symplectic structures were defined and studied from the point of view of algebraic geometry. The objects of study in this area are shifted symplectic derived stacks. As is pervasive in symplectic geometry, the study of Lagrangians is central. The notion of a Lagrangian in this context was also introduced in \cite{PTVV}. In \cite{PTVV} an interesting new structure was produced on the (homotopy) fiber product of  Lagrangians in a shifted symplectic derived stack. Namely, that given two such Lagrangians inside an $n$-shifted symplectic derived stack, their fiber product has the natural structure of an $n-1$ shifted symplectic derived stack. Our main result was a conjecture of Dominic Joyce. The proof of this result (Theorem \ref{thm:main}) uses techniques directly from the article \cite{PTVV}. Recently, an interesting preprint \cite{Cal} of Calaque appeared on categories build out of shifted symplectic derived stacks and Lagrangians in them (also Lagrangian correspondences) where some interesting new techniques were introduced. The topics in \cite{Cal} are closely related to the result which we present. We pursue these relations in separate work \cite{AB} where we look at a certain infinity category of shifted symplectic derived stacks and relate this to permutahedra inside the spaces of symplectic forms on various iterated intersections.
\section{Review of some parts the article \cite{PTVV} of Pantev, To\"{e}n, Vaqui\'{e}, and Vezzosi}\label{Section:review} As we are following closely the notation and techniques of \cite{PTVV} we refer to that work and the references contained in it for background material on derived Artin stacks $F$, tangent and cotangent complexes on them, their homotopy fiber products, and their sheaf theory, in particular the infinity categories $L_{qcoh}(F)$.  The $\infty$-category of derived stacks over $k$ for the \'{e}tale topology is denoted $\mathbf{dSt}_{k}.$ We also need the  $\infty$-categories $\mathbf{dg}_{k}^{gr}$ and $\epsilon-\mathbf{dg}_{k}$ of graded mixed complexes and of graded complexes which were described in the first section of \cite{PTVV}. We use the convention that we do not denote homotopy fiber products in any special way because all of our fiber products are homotopy fiber products of derived stacks.
Let $F$ be a derived Artin stack. Given a distinguished triangle 
\[S_1 \to S_2 \to S_3
\]
in  $L_{qcoh}(F)$, then $S_3[-1] \to S_1 \to S_2$ is also a distinguished triangle.
Let $f:F \to F'$ be a morphism of derived Artin stacks. We will use the cotangent complex $\LL_{F/k}$ of a derived Artin stack over $k$, we will simply denote this as $\LL_{F}$. More information about this can be found in HAG II \cite{TVe3} and Section 7.3 of Higher Algebra \cite{L}. There is a distinguished triangle in $L_{qcoh}(F)$ of the form 
\[f^{*}\LL_{F'}\to \LL_{F} \to \LL_{f}
\]
which defines the relative cotangent complex $\LL_{f}$. For the derived Artin stacks we will be dealing with, these complexes will be perfect and there is a dual distinguished triangle of tangent complexes 
\[ \TT_{f} \to \TT_{F} \to f^{*}\TT_{F'}.
\]
The $\infty$-category of simplicial sets is denoted $\mathbb{S}$ and these are sometimes called spaces. Given a derived Artin stack $F$, the author's of \cite{PTVV} define a space of $n$-shifted $p$-forms $\mathcal{A}^{p}(F,n) \in \mathbb{S}$ and similarly a space of $n$-shifted closed $p$-forms $\mathcal{A}^{p,cl}(F,n) \in \mathbb{S}$. We briefly outline some key steps in their definition. 
Recall that there is an $\infty$-functor 
\[NC^w :\mathbf{dSt}_{k}^{op} \longrightarrow \mathbf{dg}_{k}^{gr}
\]
defined as the composition of the $\infty$-functors
\[\mathbf{DR}:\mathbf{dSt}_{k}^{op} \longrightarrow \epsilon-\mathbf{dg}_{k}^{op}
\]
(as defined in \cite{CC} or \cite{HTT}) and a weighted negative cyclic complex functor (defined on page 15 of \cite{PTVV}) denoted
\[NC^w :\epsilon-\mathbf{dg}_{k}^{op} \longrightarrow \mathbf{dg}_{k}^{gr}.
\]

We have by Proposition 1.14 of \cite{PTVV}

\[\mathcal{A}^{p}(F,n) \cong Map_{L_{qcoh}(F)}(\mathcal{O}_{F}, \wedge^{p} \mathbb{L}_{F}[n]) \cong |\mathbf{DR}(F)[n-p](p)|
\]
and 
\[\mathcal{A}^{p,cl}(F,n)
\cong |NC^{w}(F)[n-p](p)|
\]
as $\infty$-functors $\mathbf{dSt}_{k}^{op} \to \mathbb{S}.$
There is a natural transformation
\[NC^{w}(F)[n-p](p) \to \wedge^{p}\mathbb{L}_{F}[n]
\] which gives a natural transformation of topological functors
\[\mathcal{A}^{p,cl}(F,n) \to \mathcal{A}^{p}(F,n).
\]
In Section 1.2 of \cite{PTVV} it is explained that $2$-form $\omega \in \mathcal{A}^{2}(F,n)$ gives rise \cite{PTVV} to a certain morphism
\[\Theta_{\omega}:\TT_{F} \to \LL_{F}[n].
\]
The connected components of the space $\mathcal{A}^{2}(F,n)$ for which this morphism is a quasi-isomorphism are called non-degenerate. These components form a full subspace $\mathcal{A}^{2}(F,n)^{nd}$ inside $\mathcal{A}^{2}(F,n)$. The space of $n$-shifted symplectic structures on the derived stack $F$ is defined as the homotopy pullback of $\mathcal{A}^{2}(F,n)^{nd}$ and  $\mathcal{A}^{2,cl}(F,n)$ over  $\mathcal{A}^{2}(F,n)$. Let $f:X \to F$ be a morphism of derived Artin stacks. An isotropic structure $h$ on an $n$-shifted symplectic form $\omega$ on $F$ is a path from $0$ to $f^{*}\omega$ in $X$. In Section 2.2 of \cite{PTVV} it is explained that an isotropic structure $h$ gives rise to a certain morphism 
\[\Theta_{h}: \TT_{f} \to \LL_{X}[n-1]
\]
and $h$ is said to be Lagrangian if $\Theta_{h}$ is a quasi-isomorphism.

Recall from \cite{PTVV} that given Lagrangians $f:X \to F$ and $g:Y \to F$ in an $n$-shifted symplectic derived stack $F$, corresponding to paths $h$ from $0$ to $f^{*}\omega$ and $k$ from $0$ to $g^{*}\omega$ they define an $n-1$ shifted symplectic form $R(\omega,h,k)$ on $X \times_{F} Y$ whose underlying morphism of complexes $\Theta_{R(\omega,h,k)}$ fits into a diagram with exact rows
\begin{equation}
\xymatrix@=6em{\mathbb{T}_{X \times_{F} Y}\ar[r] \ar[d]_{\Theta_{R(\omega,h,k)}}  & pr_{X}^{*} \mathbb{T}_{X} \oplus pr_{Y}^{*} \mathbb{T}_{Y} \ar[r] \ar[d]_{(\Theta_{h}, \Theta_{k})} & \pi^{*} \mathbb{T}_{F} \ar[d]_{\Theta_{\omega}}\\
\mathbb{L}_{X \times_{F} Y}[n-1]\ar[r]   & pr_{X}^{*} \mathbb{L}_{f}[n-1] \oplus pr_{Y}^{*} \mathbb{L}_{g}[n-1]\ar[r] & \pi^{*} \mathbb{L}_{F}[n]
}
\end{equation}
in which the bottom row is (a rotation of) the distinguished triangle coming from 
\begin{equation}
\xymatrix@=6em{& pr_{X}^{*} \mathbb{L}_{f}[n-1] \ar[d] \ar[r] & 0 \ar[d] \\
pr_{Y}^{*} \mathbb{L}_{g}[n-1] \ar[d] \ar[r] & \pi^{*} \LL_{F}[n] \ar[d] \ar[r]  & pr_{Y}^{*}\LL_{Y}[n] \ar[d] \\
0 \ar[r] & pr_{X}^{*}\LL_{X}[n] \ar[r] & \mathbb{L}_{X \times_{F} Y}[n]
}
\end{equation}
where the bottom square is a homotopy pushout (see \cite{TVe3}, Lemma 1.4.1.12).
\section{Multiple Derived Lagrangian Intersections}\label{MLI}

Let $(S, \omega)$ be a $n$-shifted symplectic derived stack. We will consider three Lagrangians $L,M,N$ for $S$.  Consider the homotopy cartesian diagram 
\begin{equation}\label{equation:ThreeLag}
\xymatrix{& L \times_{S} M \times_{S} N \ar[dl]_{q_N} \ar[d]_{q_M} \ar[dr]^{q_L}& \\
L \times_{S} M \ar[d]_{p_M} \ar[dr]_>>>>>>{p_L} & N \times_{S} L \ar[dr]_>>>>>>{p_L}|\hole \ar[dl]^>>>>>>{p_N}|\hole & M \times_{S} N \ar[dl]^>>>>>>{p_N} \ar[d]^{p_M} \\
L\ar[dr]_{f_L} & M \ar[d]_{f_M} & N \ar[dl]^{f_N} \\
& S &
}
\end{equation}
Note that the term $L \times_{S} M \times_{S} N$ is the homotopy limit of the rest of the diagram.
The following theorem was a conjecture of Dominic Joyce, specialized to the case of three Lagrangians. We prove the case of three Lagrangians here and then after that we will present a more general Theorem \ref{thm:main} for $m$ Lagrangians.
 \begin{thm}
 \label{thm:ThreeLag}
 Let $S$ be an $n$-shifted symplectic derived stack and let $f_L:L\to S,f_M: M\to S,f_N: N\to S$ be Lagrangian morphisms. Then 
 \[(q_N,q_L,q_M): L \times_{S} M \times_{S} N \rightarrow (L \times_{S} M) \times (M \times_{S} N) \times (N \times_{S} L)
 \] is Lagrangian with respect to the symplectic structure $q_{M}^{*}\omega_{NL}+  q_{L}^{*} \omega_{MN} +  q_{N}^{*}\omega_{LM}$, where the morphisms $q_N,q_L,q_M$ are the projection morphisms and $\omega_{NL}, \omega_{MN}, \omega_{LM}$ are the $n-1$ shifted symplectic structures provided by \cite{PTVV}. 
 \end{thm}
 {\bf Proof.}
The given Lagrangian structures for $f_L, f_M, f_N$ include isotropic structures, $h_L,h_M,h_N.$ These are 
\begin{itemize}
\item
a path $h_L$ between  $0$ and $f^{*}_{L}\omega$ in $\mathcal{A}^{2,cl}(L,n)$
\item
a path $h_M$ between  $0$ and $f^{*}_{M}\omega$ in $\mathcal{A}^{2,cl}(M,n)$
\item
a path $h_N$ between  $0$ and $f^{*}_{N}\omega$ in $\mathcal{A}^{2,cl}(N,n).$
\end{itemize}

Consider the homotopy commutativity of the lower three squares of \ref{equation:ThreeLag}. This gives 
\begin{itemize}
\item
a natural homotopy $u_{LM}$ between  $p_{M}^{*}f_{L}^{*}$ and $p_{L}^{*}f_{M}^{*}$ as maps from $\mathcal{A}^{2,cl}(S,n)$ to $\mathcal{A}^{2,cl}(L \times_{S}M,n)$ 
\item
a natural homotopy $u_{NL}$ between  $p_{L}^{*}f_{N}^{*}$ and $p_{N}^{*}f_{L}^{*}$ as maps from $\mathcal{A}^{2,cl}(S,n)$ to $\mathcal{A}^{2,cl}(N \times_{S}L,n)$ 
\item
a natural homotopy $u_{MN}$ between  $p_{N}^{*}f_{M}^{*}$ and $p_{M}^{*}f_{N}^{*}$ as maps from $\mathcal{A}^{2,cl}(S,n)$ to $\mathcal{A}^{2,cl}(M \times_{S}N,n).$
\end{itemize}
If we apply these homotopies to $\omega$ we get 
\begin{itemize}
\item a path $u_{LM}\omega$ from $p_{M}^{*}f_{L}^{*}\omega$ to $p_{L}^{*}f_{M}^{*}\omega$ in $\mathcal{A}^{2,cl}(L \times_{S}M,n)$ 
\item a path $u_{NL}\omega$ from $p_{L}^{*}f_{N}^{*}\omega$ to $p_{N}^{*}f_{L}^{*}\omega$ in $\mathcal{A}^{2,cl}(N \times_{S}L,n)$ 
\item a path $u_{MN}\omega$ from $p_{N}^{*}f_{M}^{*}\omega$ to $p_{M}^{*}f_{N}^{*}\omega$ in $\mathcal{A}^{2,cl}(M \times_{S}N,n)$ 
\end{itemize}
Similarly, the homotopy commutativity of the upper three squares of (\ref{equation:ThreeLag}) gives us
\begin{itemize}
\item
a natural homotopy $v_{L}$ between  $q_{M}^{*}p_{N}^{*}$ and $q_{N}^{*}p_{M}^{*}$ as maps from $\mathcal{A}^{2,cl}(L,n)$ to $\mathcal{A}^{2,cl}(L \times_{S}M \times_{S} N,n)$ 
\item
a natural homotopy $v_{M}$ between  $q_{N}^{*}p_{L}^{*}$ and $q_{L}^{*}p_{N}^{*}$ as maps from $\mathcal{A}^{2,cl}(M,n)$ to $\mathcal{A}^{2,cl}(L \times_{S} M\times_{S}N,n)$ 
\item
a natural homotopy $v_{N}$ between  $q_{M}^{*}p_{L}^{*}$ and $q_{L}^{*}p_{M}^{*}$ as maps from $\mathcal{A}^{2,cl}(N,n)$ to $\mathcal{A}^{2,cl}(L \times_{S} M \times_{S}N,n).$
\end{itemize}
If we apply these homotopies to the paths $h_{L}, h_{M}, h_{N}$ we get paths of paths
\begin{itemize}
\item
$k_{L}$ from $q_{N}^{*}p_{M}^{*}h_L$ to $q_{M}^{*}p_{N}^{*}h_L$
\item
$k_{M}$ from $q_{N}^{*}p_{L}^{*}h_M$ to $q_{L}^{*}p_{N}^{*}h_M$
\item
$k_{N}$ from $q_{L}^{*}p_{M}^{*}h_N$ to $q_{M}^{*}p_{L}^{*}h_N$
\end{itemize}
in the space of paths starting at the origin in $\mathcal{A}^{2,cl}(L \times_{S} M \times_{S}N,n).$ 
Consider the loops at zero
\begin{itemize}
\item
$\omega_{LM} =p_{L}^{*}h_{M}^{-1}\circ u_{LM}\omega \circ p_{M}^{*}h_{L}$ in the space $\mathcal{A}^{2,cl}(L \times_{S}M,n))$
\item
$\omega_{NL}= p_{N}^{*}h_{L}^{-1}\circ u_{NL}\omega \circ p_{L}^{*}h_{N}$ in the space $\mathcal{A}^{2,cl}(N \times_{S}L,n))$ 
\item
 $\omega_{MN}=  p_{M}^{*}h_{N}^{-1}\circ u_{MN}\omega \circ p_{N}^{*}h_{M}$ in the space $\mathcal{A}^{2,cl}(M \times_{S}N,n)).$ 
\end{itemize}
In the space $\mathcal{A}^{2,cl}(L \times_{S} M \times_{S}N,n)$ we have the figure
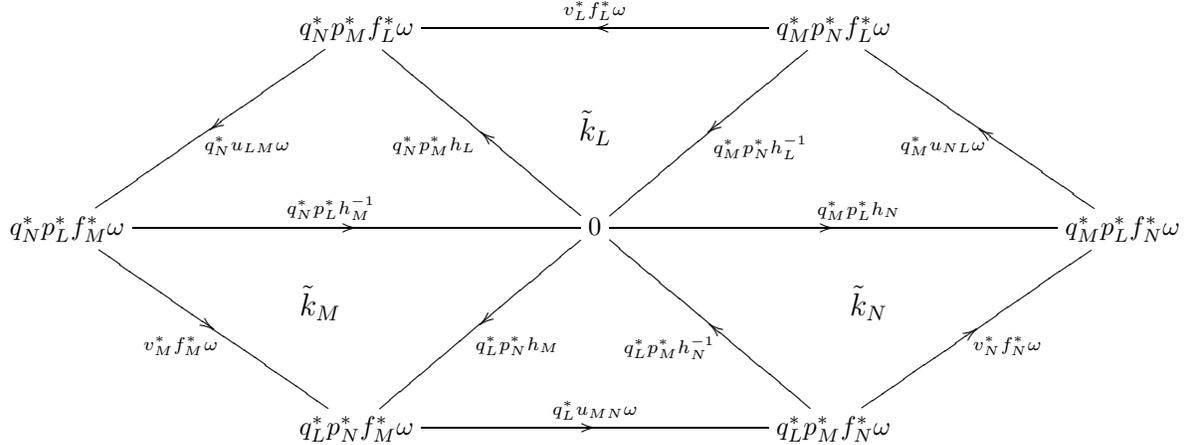
\begin{equation}\label{equation:isotropicity}
\xymatrix@=6em{
& q_{N}^{*}p_{M}^{*}f_{L}^{*}\omega  \ar@{-}|-*=0@{>}[dl]^{q_{N}^{*}u_{LM}\omega} & & q_{M}^{*}p_{N}^{*}f_{L}^{*}\omega \ar@{-}|-*=0@{>}[ll]_{v_{L}^{*}f_{L}^{*} \omega} \ar@{-}|-*=0@{>}[dl]^{q_{M}^{*}p_{N}^{*}h_{L}^{-1}}  & \\
 q_{N}^{*}p_{L}^{*}f_{M}^{*}\omega \ar@{-}|-*=0@{>}[dr]_{v_{M}^{*}f_{M}^{*}\omega} \ar@{-}|-*=0@{>}[rr]^{q_{N}^{*}p_{L}^{*}h_{M}^{-1}}
& & 0 \ar @{} [drr] ^{\mbox{\Large $\tilde{k}_{N}$}} \ar @{} [u] |{\mbox{\Large $\tilde{k}_{L}$}} \ar @{} [dll] _{\mbox{\Large $\tilde{k}_{M}$}} \ar@{-}|-*=0@{>}[rr]^{q_{M}^{*}p_{L}^{*}h_{N}} \ar@{-}|-*=0@{>}[ul]^{q_{N}^{*}p_{M}^{*}h_{L}}  \ar@{-}|-*=0@{>}[dl]^{q_{L}^{*}p_{N}^{*}h_{M}}  & & 
q_{M}^{*}p_{L}^{*}f_{N}^{*}\omega   \ar@{-}|-*=0@{>}[ul]^{q_{M}^{*}u_{NL}\omega}
\\
& q_{L}^{*}p_{N}^{*}f_{M}^{*}\omega \ar@{-}|-*=0@{>}[rr]^{q_{L}^{*}u_{MN}\omega} & & q_{L}^{*}p_{M}^{*}f_{N}^{*}\omega \ar@{-}|-*=0@{>}[ul]^{q_{L}^{*}p_{M}^{*}h_{N}^{-1}}\ar@{-}|-*=0@{>}[ur]_{v_{N}^{*}f_{N}^{*}\omega} &
}
\end{equation}
where we can recognize the concatenation 
\begin{equation}\label{equation:concat}q_{M}^{*}\omega_{NL} \circ q_{L}^{*} \omega_{MN}\circ q_{N}^{*}\omega_{LM}
\end{equation}
of loops in $\mathcal{A}^{2,cl}(L \times_{S} M \times_{S}N,n).$ In this figure, using $k_{L}$, $k_{M}$, and $k_{N}$ we filled in three triangles with $2$-simplices which have been denoted $\tilde{k}_{L}$, $\tilde{k}_{M}$, and $\tilde{k}_{N}$.
 
The fact that (\ref{equation:ThreeLag}) is homotopy cartesian implies that the path around the outside of (\ref{equation:isotropicity}) is the boundary of a filled in hexagon in the space $\mathcal{A}^{2,cl}(L \times_{S} M \times_{S}N,n).$ Therefore we have produced a homotopy from the constant loop at $0$ to the the concatenation of loops (\ref{equation:concat}). This can be interpreted as a homotopy from the constant loop at $0$ to the loop $q_{M}^{*}\omega_{NL} \circ q_{L}^{*} \omega_{MN}\circ q_{N}^{*}\omega_{LM}$ at $0$ in the space $\mathcal{A}^{2,cl}(L \times_{S} M \times_{S}N,n).$  Because concatenation and sum agree up to homotopy this can be interpreted as a homotopy from the constant loop at $0$ to the loop $q_{M}^{*}\omega_{NL}+  q_{L}^{*} \omega_{MN} +  q_{N}^{*}\omega_{LM}.$

This gives us a path $h$ from $0$ to the corresponding element of the space $\mathcal{A}^{2,cl}(L \times_{S} M \times_{S}N,n-1)$. But this element is precisely the pullback of the product symplectic structure on $(L \times_{S} M) \times (M \times_{S} N) \times (L \times_{S} N)$ to $L \times_{S} M \times_{S} N$ given by the construction from \cite{PTVV} applied to each of the three components and described in Section \ref{Section:review}.
Therefore, we have produced an isotropic structure $h$ on the morphism 
\[ L \times_{S} M \times_{S} N
 \to (L \times_{S} M) \times (M \times_{S} N) \times (N \times_{S} L).\]

\noindent
Let 
\[r_{L}: L \times_{S} M \times_{S}N \to L
\]
\[r_{M}: L \times_{S} M \times_{S}N \to M
\]
\[r_{N}: L \times_{S} M \times_{S}N \to N
\]
and 
\[r: L \times_{S} M \times_{S}N \to S
\]
be the projection morphisms.
 In order to show that the isotropic structure $h$ is Lagrangian, we need to consider a pair of $4$ term rows which correspond to diagram (\ref{equation:ThreeLag}). By definition of the isotropic structure $h$ we have a commutative diagram where the top row is the homotopy limit of the rest of the diagram:
\begin{equation}
\xymatrix{ \mathbb{T}_{L \times_{S} M \times_{S} N} \ar[d] \ar[rrrr]^{\Theta_{h}} &&&& \mathbb{L}_{(q_N,q_M,q_L)}[n-2] \ar[d] \\
q^{*}_N \mathbb{T}_{L \times_{S} M} \oplus q^{*}_M \mathbb{T}_{N \times_{S} L} \oplus q^{*}_L \mathbb{T}_{M \times_{S} N} \ar[d] \ar[rrrr]^{\Theta_{q_{M}^{*}\omega_{NL}+  q_{L}^{*} \omega_{MN} +  q_{N}^{*}\omega_{LM}}} &&&& (q^{*}_N \mathbb{L}_{L \times_{S} M} \oplus q^{*}_M \mathbb{L}_{N \times_{S} L} \oplus q^{*}_L \mathbb{L}_{M \times_{S} N})[n-1] \ar[d] \\
 r_L^*\mathbb{T}_{L} \oplus r_M^*\mathbb{T}_{M} \oplus r_N^*\mathbb{T}_{N} \ar[d] \ar[rrrr]^{(\Theta_{h_L},\Theta_{h_M},\Theta_{h_N}) } &&&& (r_L^* \mathbb{L}_{f_L} \oplus r_M^* \mathbb{L}_{f_M} \oplus r_N^*\mathbb{L}_{f_N})[n-1] \ar[d] \\
r^{*} \mathbb{T}_S \ar[rrrr]^{\Theta_\omega} &&&&  r^{*} \mathbb{L}_S [n]
}
\end{equation}

 where
 \[\Theta_{q_{M}^{*}\omega_{NL}+  q_{L}^{*} \omega_{MN} +  q_{N}^{*}\omega_{LM}}= (\Theta_{R(\omega_{LM},h_L,h_M)}, \Theta_{R(\omega_{NL},h_N,h_L)}, \Theta_{R(\omega_{MN},h_M,h_N)}).
 \]
 Due to the fact that the lower three horizontal arrows are quasi-isomorphisms, the top horizontal arrow is as well. Therefore, the isotropic structure $h$ which we constructed endows the morphism $(q_{N}, q_{L}, q_{M})$ with a Lagrangian structure.
\ \hfill $\Box$
\bigskip

\noindent  
Similarly, we have the following Theorem (conjectured by Dominic Joyce)
\begin{thm}\label{thm:main}
Let $m$ be an integer greater than equal to $2$. Let $S$ be an $n$-shifted symplectic derived stack and let $L_i\to S$ be Lagrangian morphisms for $i=1, \dots, m$. Then the morphism
  \[L_1 \times_{S} \cdots \times_{S} L_m \to (L_1 \times_{S} L_2) \times (L_2 \times_{S} L_3) \times \cdots \times (L_{m-1} \times_{S} L_m) \times (L_m \times_{S} L_1)
  \]
  is Lagrangian with respect to the product $n-1$-shifted symplectic structure on the right hand side provided by \cite{PTVV}.
 \end{thm}
 {\bf Proof.}
When $m \geq 3$, the proof follows precisely along the lines of the proof of Theorem \ref{thm:ThreeLag}. The only difference in the proof is notational, for instance, the analogue of the figure in (\ref{equation:isotropicity}) now must be divided into $2m$ triangles. When $m=2$ this is just the diagonal morphism to the product of $L_1 \times_{S} L_2$ and its opposite. It is easy to see that this is Lagrangian.
 \ \hfill $\Box$

\section{Examples}
In this section we mention some explicit examples of Lagrangians in  shifted symplectic derived stack. The idea of these Lagrangians goes back to Tyurin (see for example \cite{Ty}) and was recently explained in \cite{Cal} by Calaque. Tyurin's construction is in fact a holomorphic (or algebraic) version of an older construction by Casson relating to moduli of local systems. This construction became popular because of its role in Donaldson-Thomas theory (\cite{Th}, \cite{ThThesis}). Derived stacks are expected to play an important role in a deeper understanding of Donaldson--Thomas theory, see for example \cite{BBJ, BBDJS, BBBJ, J, To, Hua, KoSo1, KoSo2, KL} and the references therein.

One type of examples come from Calabi-Yau $3$-folds which degenerate into a pair of Fano $3$-folds intersecting in a $K3$ surface. Sometimes, one passes to the stack of expanded degenerations \cite{LW}. These degenerations provide other examples of the same type. One would like to extract invariants from the derived moduli stacks of complexes of sheaves on the generic (smooth) fiber and argue that these invariants are constant in this degenerating family. Therefore, these invariants can be extracted from some derived moduli stacks of complexes of sheaves on the singular (special) fiber. One would like to relate these invariants to the so-called relative invariants of the Fano 3-folds relative to the $K3$ surface. For more details see the article \cite{LW} by J. Li and B. Wu and the references therein.

The key feature of this example (discussed in \cite{Ty} and \cite{Cal}) is that the inclusion of the $K3$ surface inside the Fano (or similar manifold) induces a Lagrangian strucutre on the derived stack of perfect complexes. One can ask for which other types of morphisms of derived Artin stacks induce Lagrangian morphisms on their derived stacks of perfect complexes. This article can be used to construct many other examples of this phenomenon. For example,
\[(F_{1}\coprod_{Z}F_{2}) \coprod (F_{2}\coprod_{Z}F_{3}) \coprod (F_{3}\coprod_{Z}F_{1} ) 
\] 
mapping to the homotopy colimit of
\begin{equation}\label{equation:ThreePreLag}
\xymatrix{
F_{1} \coprod_{Z} F_{2}  &  F_{1} \coprod_{Z} F_{3} & F_{2} \coprod_{Z} F_{3}   \\
F_{1}\ar[u] \ar[ur]|\hole & F_{2}\ar[ur] \ar[ul] & F_{3} \ar[ul]|\hole \ar[u] \\
& Z \ar[ur] \ar[ul] \ar[u] &
}
\end{equation}
gives a Lagrangian structure after passing to derived stacks of perfect complexes whenever the morphisms $Z \to F_i$ have this same property for $i=1,2,3$. In otherwords, we only need to assume that there is a shifted symplectic structure on the stack of perfect complexes on $Z$ and that the morphisms $Z \to F_i$ induce Lagrangian morphisms on derived stacks of perfect complexes.

\end{document}